\documentclass{fakeetds}         \usepackage{amsfonts} \usepackage{amssymb} \usepackage{latexsym} \usepackage{graphicx}  \usepackage{amsbsy} \usepackage{psfrag}    
   \newcommand{\D}{\displaystyle}
\newcommand{\restr}[1]{\raisebox{-0.3em}{$\lb|_{#1}\rb.$}} 
\newcommand{\ignore}[1]{}    
\newcommand{\breath}{\medskip} 
\newtheorem{thm}{Theorem}

\newcounter{claimcount}[thm]  
\newtheorem{prop}[thm]{Proposition} 

\newtheorem{lemma}[thm]{Lemma} 

\newtheorem{cor}[thm]{Corollary}

\newcommand{\dfn}{\sf\em} 
\newcommand{\Theorem}[2]{\begin{thm}{\sf #1}  #2 \end{thm}}
\newcommand{\Proposition}[2]{\begin{prop}{\sf #1}  #2 \end{prop}}
\newcommand{\Lemma}[2]{\begin{lemma}{\sf #1}  #2 \end{lemma}}
\newcommand{\Corollary}[2]{\begin{cor}{\sf #1}  #2 \end{cor}} 
\newcommand{\thmfont}[1]{{\sl #1}}    
\newcommand{\example}[1]{        \refstepcounter{thm}                     \begin{list}{} 			{\setlength{\leftmargin}{0em} 			\setlength{\rightmargin}{0em}}        \item {\sc Example \thethm:} #1                   \hfill$\diamondsuit$  \end{list}   			}     
\newcommand{\bthmlist}{ \begin{list}{{\bf(\alph{enumi})}} {\usecounter{enumi} \setlength{\leftmargin}{1em} \setlength{\itemsep}{0.2em} \setlength{\topsep}{0.2em} \setlength{\itemindent}{0em} \setlength{\parsep}{0em} \setlength{\rightmargin}{0em}} } 
\newcommand{\ethmlist}{\end{list}}    
\newcommand{\Claim}[1]{\refstepcounter{claimcount}                \noindent {\sc Claim \theclaimcount: \ }\thmfont{ #1}} 
\newcommand{\bprf}[1][Proof:]{\begin{list}{} 			{\setlength{\leftmargin}{0.7em} 			\setlength{\rightmargin}{0em} 			\setlength{\listparindent}{1em}}                         \item {\em \hspace{-1em}  #1  }} 
\newcommand{\eprf}{\end{list}} 
\newcommand{\bthmprf}{\bprf}
\newcommand{\bclaimprf}{\bprf}
\newcommand{\ethmprf}{ \hfill$\Box$  \eprf  \breath  } 
\newcommand{\eclaimprf}{ \hfill $\Diamond$~{\scriptsize {\tt Claim~\theclaimcount}}\eprf}  
\newcommand{\QED}{\hfill\ensuremath{\Box}}
\newcommand{\qed}{\QED}     
\newcommand{\beq}{\begin{eqnarray*}}
\newcommand{\eeq}{\end{eqnarray*}} 
\newcommand{\beqn}{ \begin{equation} }
\newcommand{\eeqn}{ \end{equation} }
\newcommand{\bdesc}{\begin{description}}
\newcommand{\edesc}{\end{description}}   
\newcommand{\If}{\mbox{\ if \ }}  
\newcommand{\dB}{{\mathbb{B}}}
\newcommand{\dC}{{\mathbb{C}}}
\newcommand{\dF}{{\mathbb{F}}}
\newcommand{\dH}{{\mathbb{H}}}
\newcommand{\dK}{{\mathbb{K}}}
\newcommand{\dM}{{\mathbb{M}}}
\newcommand{\dN}{{\mathbb{N}}}
\newcommand{\dT}{{\mathbb{T}}}
\newcommand{\dU}{{\mathbb{U}}}
\newcommand{\dV}{{\mathbb{V}}}
\newcommand{\dW}{{\mathbb{W}}}
\newcommand{\dZ}{{\mathbb{Z}}}       
\newcommand{\bB}{{\mathbf{ B}}}
\newcommand{\bC}{{\mathbf{ C}}}
\newcommand{\bG}{{\mathbf{ G}}}
\newcommand{\bS}{{\mathbf{ S}}}
\newcommand{\bU}{{\mathbf{ U}}}
\newcommand{\ba}{{\mathbf{ a}}}
\newcommand{\bb}{{\mathbf{ b}}}
\newcommand{\bc}{{\mathbf{ c}}}
\newcommand{\bg}{{\mathbf{ g}}}
\newcommand{\bh}{{\mathbf{ h}}}
\newcommand{\bk}{{\mathbf{ k}}}
\newcommand{\bs}{{\mathbf{ s}}}
\newcommand{\bu}{{\mathbf{ u}}}
\newcommand{\bchi}{{\boldsymbol{\chi }}}
\newcommand{\sA}{{\mathcal{ A}}}
\newcommand{\sB}{{\mathcal{ B}}}
\newcommand{\sE}{{\mathcal{ E}}}
\newcommand{\sF}{{\mathcal{ F}}}
\newcommand{\sI}{{\mathcal{ I}}}
\newcommand{\sR}{{\mathcal{ R}}}
\newcommand{\gam }{\gamma}
\newcommand{\sig }{\sigma} 
\newcommand{\Gam }{\Gamma}
\newcommand{\hbG}{{\widehat{\mathbf{ G}}}}
\newcommand{\hbS}{{\widehat{\mathbf{ S}}}}
\newcommand{\hsA}{{\widehat{\mathcal{ A}}}}
\newcommand{\hmu }{{\widehat{\mu}}}
\newcommand{\hnu }{{\widehat{\nu}}}
\newcommand{\heta}{{\widehat{\eta }}} 
\newcommand{\fb}{{\mathsf{ b}}}
\newcommand{\fh}{{\mathsf h}} 
\newcommand{\fk}{{\mathsf{ k}}}
\newcommand{\fm}{{\mathsf{ m}}}
\newcommand{\fn}{{\mathsf{ n}}}
\newcommand{\fp}{{\mathsf{ p}}}
\newcommand{\fu}{{\mathsf{ u}}}
\newcommand{\fv}{{\mathsf{ v}}}
\newcommand{\tla}{{\widetilde{a}}}
\newcommand{\tlr}{{\widetilde{r}}}
\newcommand{\tlbZ}{{\widetilde{\mathbf{ Z}}}} 
\newcommand{\tlsA}{{\widetilde{\mathcal{ A}}}}
\newcommand{\tlsR}{{\widetilde{\mathcal{ R}}}}
\newcommand{\lb}{\left}
\newcommand{\rb}{\right} 
\newcommand{\maketall}{\rule[-0.5em]{0em}{1em}}        
\newcommand{\map}{{\longrightarrow}}
\newcommand{\goto}{{\rightarrow}}
\newcommand{\into}{{\map}}
\newcommand{\seilpmi}{{\Longleftarrow}}
\newcommand{\statement}[1]{\lb(  \maketall       \begin{minipage}{40em}       \begin{tabbing}         #1        \end{tabbing}      \end{minipage}  \rb)}     
\newcommand{\oo}{{\infty}}        
\newcommand{\X}{\times}
\newcommand{\x}{\X}
\newcommand{\dirsum}{\oplus}
\newcommand{\Dirsum}{\bigoplus} 
\newcommand{\intsct}{\cap}
\newcommand{\Intsct}{\bigcap}
\newcommand{\set}[2]{{\left\{ #1 \; ; \; #2 \right\} }} 
\newcommand{\supp}[1]{{\sf supp}\lb(#1\rb)}      
\newcommand{\inn}[1]{{\left\langle #1 \right\rangle }}       
\newcommand{\Id}[1]{{\mathbf{ Id}_{{#1}}}}
\newcommand{\chr}[1]{{{\boldsymbol{1}}_{{#1}}}} 
\newcommand{\choice}[1]{{\lb\{ \begin{array}{rcl}                                 #1                                \end{array}  \rb.  }}                     %
\newcommand{\eeequals}[1]{\raisebox{-0.9ex}{$\overline{\overline{{\scriptscriptstyle{\mathrm{#1}}}}}$}} 
\newcommand{\leeeq}[1]{\raisebox{-1ex}{${{\D\leq} \atop {\scriptscriptstyle{\mathrm{#1}}}}$}} 
\newcommand{\grt}[1]{\raisebox{-1ex}{${{\D>} \atop {\scriptscriptstyle{\mathrm{#1}}}}$}}  
\newcommand{\geeeq}[1]{\raisebox{-1ex}{${{\geq} \atop {\scriptscriptstyle{\mathrm{#1}}}}$}} 
\newcommand{\Fix}[1]{{\sf Fix}\lb[#1\rb]}
\newcommand{\shift}[1]{\sig^{#1}}    
\newcommand{\dmu}{{ \;\; d\mu}}  
\newcommand{\End}[2][]{\mathsf{End}_{#1} \lb(#2\rb)}
\newcommand{\Natur}{\dN}
\newcommand{\Zahl}{\dZ}
\newcommand{\Zahlmod}[1]{{\Zahl_{/#1}}}
\newcommand{\Cplx}{\dC}
\newcommand{\CC}[1]{{\lb[ #1 \rb]}}
\newcommand{\CO}[1]{{\lb[ #1 \rb)}}

\renewcommand{\implies}{\ensuremath{\Longrightarrow}}
\renewcommand{\And}{\mbox{\ and \ }} 
\newcommand{\bvphi}{\overline{\varphi}}
\newcommand{\DE}{{D\!+\!E}} 
\newcommand{\hAM}{\widehat{\AM}} 
\newcommand{\AM}{\sA^\dM}
\newcommand{\AjM}{\sA_j^\dM}
\newcommand{\AiM}{\sA_i^\dM}
\newcommand{\AH}{\sA^\dH}
\newcommand{\AB}{\sA^\dB}
\newcommand{\Nh}{\dH}
\newcommand{\ZD}[1][D]{\Zahl^{#1}}
\newcommand{\NE}[1][E]{\Natur^{#1}}
\newcommand{\AZD}[1][D]{\sA^{\ZD[#1]}}
\newcommand{\AZ}{\sA^{\Zahl}}
\newcommand{\AN}{\sA^{\Natur}}   
\newcommand{\cyl}[1]{\inn{#1}} 
\newcommand{\htop}{h_{\mathrm{top}}}
\newcommand{\bzero}{\boldsymbol{0}} 
\newcommand{\Meas}{\mathfrak{M\scriptscriptstyle{\!e\!a\!s}}} 
\newcommand{\Fpk}{\dF_{\!p^k}} 
\newcommand{\CA}[1][\AM]{\mathsf{C\!A}(#1)}
\newcommand{\RLCA}[1][\AM]{\sR\mbox{-}\mathsf{LCA}\lb(#1\rb)}
\newcommand{\RjLCA}[1][\sA_j^\dM]{\sR_j\mbox{-}\mathsf{LCA}\lb(#1\rb)}
\newcommand{\ZpLCA}[1][\AM]{\Zahlmod{p}\mbox{-}\mathsf{LCA}\lb(#1\rb)}
\newcommand{\ZmLCA}[1][\AM]{\Zahlmod{m}\mbox{-}\mathsf{LCA}\lb(#1\rb)}
\newcommand{\ECA}[1][\AM,\star]{\mathsf{ECA}\lb(#1\rb)} 
\newcommand{\subseteeeq}[1]{\raisebox{-1.3ex}{$\stackrel{\D\subseteq}{\scriptscriptstyle{\mathrm{#1}}}$}} 
\newcommand{\iiin}[1]{\raisebox{-1.3ex}{$\stackrel{\D\in}{\scriptscriptstyle{\mathrm{#1}}}$}} 

\begin{document}

\ETDS{1}{40}{26}{2006} 
\title{Module Shifts and Measure Rigidity in Linear Cellular Automata}
\runningheads{M. Pivato}{Module shifts and Measure Rigidity in LCA}
\author{Marcus Pivato}
\address{Dept. of Mathematics, Trent University \\
1600 West Bank Drive, Peterborough, Ontario, K9J 7B8, Canada.}

\email{marcuspivato@trentu.ca}
\recd{10 July 2007}

\begin{abstract}
Suppose $\sR$ is a finite commutative ring of prime characteristic, $\sA$
is a finite $\sR$-module, $\dM:=\ZD\x\NE$, 
and $\Phi$ is an $\sR$-linear cellular
automaton on $\AM$.  If $\mu$ is a $\Phi$-invariant measure which is
multiply $\shift{}$-mixing in a certain way, then we show that $\mu$
must be the Haar measure on a coset of some submodule shift of $\AM$.
Under certain conditions, this means $\mu$ must be the uniform
Bernoulli measure on $\AM$.
\end{abstract}

  Let $\sA$ be a finite set.  Let $\dM:=\ZD\x\NE$ be a
$(\DE)$-dimensional lattice, for some $D,E\in\Natur$, and let $\AM$
denote the set of all functions $\ba:\dM\into\sA$, which we regard as
$\dM$-indexed {\dfn configurations} of elements in $\sA$.  We write
such a configuration as $\ba=[a_\fm]_{\fm\in\dM}$, where $a_\fm\in\sA$
for all $\fm\in\dM$.  Treat $\sA$ as a discrete topological space;
then $\AM$ is a Cantor space ---i.e. it is compact, perfect, totally
disconnected, and metrizable.

 If $\ba\in\AM$ and $\dU\subset\dM$, then we define
$\ba_\dU\in\sA^\dU$ by $\ba_\dU:=[a_\fu]_{\fu\in\dU}$.  If
$\fm\in\dM$, then strictly speaking, $\ba_{\fm+\dU}\in\sA^{\fm+\dU}$;
however, it will often be convenient to `abuse notation' and treat
$\ba_{\fm+\dU}$ as an element of $\sA^{\dU}$ in the obvious way.  Let
$\Nh\subset\dM$ be some finite subset, and let $\phi:\sA^\Nh\into\sA$
be a function (called a {\dfn local rule}).  The {\dfn cellular
automaton} (CA) determined by $\phi$ is the function
$\Phi:\AM\into\AM$ defined by $\Phi(\ba)_{\fm} = \phi(\ba_{\fm+\Nh})$
for all $\ba\in\AM$ and $\fm\in\dM$. We refer to $\Nh$ as the {\dfn
neighbourhood} of $\Phi$. 

  We will prove a new `measure rigidity' result for linear CA: if
$\Phi$ is a linear CA and $\mu$ is a $\Phi$-invariant measure which is
multiply $\shift{}$-mixing in a certain way, then $\mu$ must be the
Haar measure on a coset of some submodule shift of $\AM$.
In particular, if $\AM$ admits
no proper mixing subgroup shifts (e.g $D=1$ and $\sA=\Zahlmod{p}$, for
$p$ prime), then $\mu$ must be the uniform measure on $\AM$.  This
result is complementary to previous rigidity results of
\cite{Schmidt95b,HostMaassMartinez,PivatoPerm,Einsiedler05,Sablik}.

\breath

{\em Terminology \& Notation.} \  Throughout, lowercase bold-faced letters ($\ba,\bb,\bc,\ldots$) denote
elements of $\AM$, and Roman letters ($a,b,c,\ldots$) are elements of
$\sA$ or ordinary numbers.  Lower-case sans-serif
($\ldots,\fm,\fn,\fp$) are elements of $\dM$, and upper-case hollow
font ($\dU,\dV,\dW,\ldots$) are subsets of $\dM$.
 For any $\fv\in\dM$, let $\shift{\fv}:\AM\into\AM$ be the {\dfn shift
map} defined by $\shift{\fv}(\ba)_{\fm} = a_{\fm+\fv}$ for all
$\ba\in\AM$ and $\fm\in\dM$.  Let $\CA$ denote the set of cellular
automata on $\AM$; then $\CA$ is also the set of continuous
transformations of $\AM$ which commute with all shifts \cite[Theorem
3.4]{Hedlund}.  A {\dfn subshift} is a closed subset $\bS\subseteq\AM$
which is invariant under all shifts.  Let
$\CA[\bS]:=\set{\Phi\in\CA}{\Phi(\bS)\subseteq\bS}$.  

  If $\ba\in\AM$, and $\dK\subset\dM$, recall that
$\ba_\dK:=[a_\fk]_{\fk\in\dK}\in\sA^\dK$.  If $\bS\subseteq\AM$ is a
subshift, let $\bS_\dK:=\set{\bs_\dK}{\bs\in\bS}\subseteq\sA^\dK$.  If
$\bk\in\sA^\dK$, then let $\cyl{\bk}:=\set{\ba\in\AM}{\ba_\dK=\bk}$ be
the {\dfn cylinder set} defined by $\bk$.  The topology (and hence,
the Borel sigma-algebra) of $\AM$ is generated by the collection of
all such cylinder sets for all finite $\dK\subset\dM$.  Let
$\Meas(\AM)$ [resp. $\Meas(\bS)$] 
be the set of Borel probability measures on $\AM$ [resp. $\bS$], and
let $\Meas(\AM,\shift{})$ [resp. $\Meas(\bS,\shift{})$]
 be the shift-invariant measures on $\AM$ [resp. $\bS$].
If $\Phi\in\CA[\bS]$, let $\Meas(\bS,\Phi)$ be the $\Phi$-invariant
measures on $\bS$.

\breath {\em Linear CA.} \ 
Let $(\sR,+,\cdot)$ be a finite
ring with unity $1_\sR$, and let $(\sA,+,\cdot)$ be a finite
$\sR$-module.  If $\Phi\in\CA$, then $\Phi$ is an {\dfn
$\sR$-linear} CA ($\sR$-LCA) if the local rule $\phi:\AH\into\sA$ has the form
 \beqn
\label{linear.CA.local.rule}
  \phi(\ba_\Nh) \quad:=\quad \sum_{\fh\in\Nh} \varphi_\fh a_\fh,
\qquad \forall \ \ba_\Nh \in\sA^\Nh,
\eeqn
where $\varphi_\fh\in\sR\setminus\{0\}$ for each $\fh\in\Nh$.  
Let $\RLCA$ be the set of all $\sR$-linear CA on $\AM$.

\example{\label{X:LCA} If $(\sA,+)$ is a finite abelian group,
 and $\AM$ is treated as
a Cartesian product and endowed with componentwise addition,
then $(\AM,+)$ is a compact abelian group.  If $\Phi\in\CA$, then
$\Phi$ is an {\dfn endomorphic cellular automaton} (ECA) if $\Phi$ is
also a group homomorphism of $(\AM,+)$.  Let $\sE$ be the
(noncommutative) ring of all group endomorphisms of $(\sA,+)$.  Then
$\sA$ is an $\sE$-module, and any ECA on $\AM$ is an $\sE$-linear CA.  (In the
literature, these are often just called {\em linear} CA.)

(b) In particular, let  $m\in\Natur$ and let $\sA=(\Zahlmod{m},+)$,
with addition modulo $m$.
Then $\sE=(\Zahlmod{m},+,\cdot)$ [with multiplication modulo $m$], so
eqn.(\ref{linear.CA.local.rule}) becomes 
$\phi(\ba_\Nh)  =  \sum_{\fh\in\Nh} \varphi_\fh a_\fh \ \bmod{m}$,
where $\varphi_\fh\in\CO{1...m}$ for each $\fh\in\Nh$.   

(c) Let $k\in\Natur$, and let $\sR:=\Fpk$ be the unique finite
field of order $p^k$ [in particular, if $k=1$ then $\Fpk=\Zahlmod{p}$ as in
Example (b)].  Let $\sA$ be any finite-dimensional
$\Fpk$-vector space (e.g. $\sA:=(\Fpk)^m$, for some $m\in\Natur$);
  then $\AM$ is an (infinite-dimensional)
$\Fpk$-vector space, and $\Phi$ is an $\Fpk$-LCA iff $\Phi$
is a linear endomorphism of $\AM$.

(d) Let $\sA$ be a finite-dimensional $\Fpk$-vector space, and let
$\End{\sA}$ be the (noncommutative) ring of all $\Fpk$-linear
endomorphisms of $\sA$.  Suppose
$\{\varphi_\fh\}_{\fh\in\dH}\subset\End{\sA}$ is a collection of
endomorphisms which commute with one another, and let $\sR$ be the
subring of $\End{\sA}$ generated by $\{\varphi_\fh\}_{\fh\in\dH}$;
then $\sR$ is a commutative ring, $\sA$ is an $\sR$-module, and if
$\Phi$ is as in eqn.(\ref{linear.CA.local.rule}), then $\Phi$ is an
$\sR$-LCA.}

Example \ref{X:LCA}(a) is `universal' in the following sense: any
$\sR$-module is also an abelian group, and any $\sR$-linear CA is
automatically an  ECA.  However, in a general ECA, the
coefficients $\{\varphi_\fh\}_{\fh\in\dH}\subset\sE$ do not
commute, because the endomorphism ring $\sE$ is not commutative
unless $\sA=\Zahlmod{m}$, as in Example \ref{X:LCA}(b).  If
$\sR$ is commutative, then $\sR$-LCA are much easier to analyze
than general ECA.

  If $r\in\sR$, then the {\dfn characteristic} of $r$ is the smallest
$m\in\Natur$ such that $m\cdot r=0$ (or it is $0$ if there is no such
$m$).  The {\dfn characteristic} of $\sR$ is the characteristic of the
unity element $1_\sR$.  (For example, $\sR=\Zahlmod{m}$ has
characteristic $m$.)  If $\sR$ has characteristic $m$, then $m\cdot
r=0$ for all $r\in\sR$; hence the characteristic of $r$ divides $m$.
We will be mainly interested in the case when $\sR$ is a {\em
commutative} ring of prime characteristic, as in Examples
\ref{X:LCA}(c,d).

\breath

{\em Subgroup shifts and submodule shifts.} \ 
 Suppose $(\sA,+)$ is a finite abelian group, so that
$(\AM,+)$ is compact abelian.  A {\dfn subgroup shift}
is a closed, shift-invariant subgroup $\bG\subset \AM$ (i.e. $\bG$ is
both a subshift and a subgroup); see
\cite{Kitchens1,Kitchens2,KitchensSchmidt89,KitchensSchmidt92,Schmidt95}.
If $\sA$ is an $\sR$-module, then $\AM$ is also an $\sR$-module under
componentwise $\sR$-multiplication.  An {\dfn $\sR$-submodule shift}
is a subgroup shift which is also an $\sR$-submodule. 
For example, if $\Phi\in\RLCA$, then $\Phi(\AM)$ and 
 $\ker(\Phi):=\Phi^{-1}\{\bzero\}$ are
submodule shifts (here $\bzero\in\AM$ is the constant zero element).
Also $\Fix{\Phi}:=\set{\ba\in\AM}{\Phi(\ba)=\ba}$
is a submodule shift (because $\Fix{\Phi}=\ker(\Phi-\Id{})$). 
If $\sA=\sR=\Zahlmod{m}$, then every subgroup shift is a submodule shift,
and vice versa.  However, in general the $\sR$-submodule shifts form a
more restricted class.

To study the ergodic theory of $\sR$-LCA, it is first necessary to
characterize their invariant measures. 
If $\bG\subseteq\AM$ is a subgroup shift, then the
{\dfn Haar measure} of
$\bG$ is the unique $\eta_\bG\in\Meas(\bG)$ which is
invariant under translation by all elements of $\bG$.  That is, if
$\bg\in\bG$, and $\bU\subset\bG$ is any measurable subset, and
$\bU+\bg:=\set{\bu+\bg}{\bu\in\bU}$, then
$\eta_\bG[\bU+\bg]=\eta_\bG[\bU]$.  In particular, if $\bG=\AM$, then
$\eta_\bG$ is just the uniform Bernoulli measure on $\AM$. 
Let $\ECA[\bG]:=
\{\Phi\in\CA \ ; $ $\Phi$ is an ECA and $\Phi(\bG)\subseteq\bG\}$.

\Proposition{\label{subgroup.haar.invariant}}
{
  Let $(\sA,+)$ be a finite abelian group, let $\bG\subseteq\AM$ be
a subgroup shift, and let $\Phi\in\ECA[\bG]$.  Then
  $\statement{$\Phi(\eta_\bG)=\eta_\bG$} \iff \statement{$\Phi(\bG)=\bG$}$. 
}
\bthmprf
{\bf(a)}  ``$\implies$'' is because $\bG=\supp{\eta_\bG}$.
To see  ``$\seilpmi$'', note that $\Phi(\eta_\bG)\in\Meas(\bG)$.
Thus, it suffices to show that $\Phi(\eta_\bG)$ is invariant under all
$\bG$-translations.  Let $\bg\in\bG$, and let $\tau^\bg:\bG\into\bG$
be the translation map [i.e. $\tau^\bg(\bh):=\bg+\bh$].
Find $\bh\in\bG$ such
that $\Phi(\bh)=\bg$ (this $\bh$ exists because $\Phi(\bG)=\bG$).
Then $\tau^\bg \lb[\Phi(\eta_\bG)\rb] \ \eeequals{(*)} \
 \Phi\lb[\tau^\bh(\eta_\bG)\rb]\ \eeequals{(\dagger)} \
 \Phi(\eta_\bG)$, where $(*)$ is because 
$\tau^\bg \circ\Phi = \Phi\circ\tau^\bh$,
and $(\dagger)$ is because $\eta_\bG$ is the Haar measure.
This holds for all $\bg\in\bG$.  But $\eta_\bG$ is the unique
probability measure on $\bG$
 such that $\tau^\bg(\eta_\bG)=\eta_\bG$ for all $\bg\in\bG$;
thus $\Phi(\eta_\bG)=\eta_\bG$.\ethmprf

  Some ECA exhibit a great deal of {\em measure rigidity}, meaning that
the Haar measures of $\Phi$-invariant subgroup shifts are 
the {\em only} $\Phi$-invariant measures satisfying certain
`nondegeneracy' conditions.  For example, Host, Maass and Mart\'inez
showed that, if $p$ is prime and $\sA=\Zahlmod{p}$, and $\Phi$ is a
nearest neighbour $\Zahlmod{p}$-LCA on $\AZ$, then the only
positive-entropy, $\shift{}$-ergodic, $\Phi$-invariant measure is 
 the Haar measure on $\AZ$ ---i.e. the uniform Bernoulli measure 
\cite[Thm.12]{HostMaassMartinez}.  This result was vastly generalized
by Sablik, who showed that, if $\bG\subseteq\AZ$ is any subgroup shift
and $\Phi\in\ECA[\bG]$, then the only positive-entropy,
$\Phi$-invariant measure satisfying certain  ergodicity
conditions is the Haar measure on $\bG$ \cite[Thm. 3.3 and
3.4]{Sablik}. (Actually, Sablik's result is even more general, since it
allows any abelian group shift structure on $\AZ$.)  See also
\cite{PivatoPerm} for similar results concerning ECA in the
case when $\AZ$ is a {\em nonabelian} group shift, as well as {\em
multiplicative} CA (in the case when $(\sA,\cdot)$ is a nonabelian
group).

  All of these results are for {\em one}-dimensional ECA.  Einsiedler
 \cite[Corollary 2.3]{Einsiedler05} has a similar rigidity result for
 automorphic $\ZD$-actions on compact abelian groups (e.g. the $\ZD$-shift
 action on a subgroup shift $\bG\subseteq\AZD$).  This theorem can be
 easily translated into equivalent rigidity results for ECA in
 $\AZD[D-1]$.  Like \cite{HostMaassMartinez} and
 \cite{Sablik}, Einsiedler requires both an entropy condition and
 fairly strong ergodicity hypotheses.

\breath

  Let $\sR$ be a commutative ring of characteristic $p$.  We will
prove a measure rigidity result for multidimensional $\sR$-LCA whose
only requirement on the measure is a limited form of multiple mixing.
Our result is philosophically similar to the rigidity results in
\cite{Schmidt95b} or \cite[\S29]{Schmidt95}, but it is applicable to
much larger class of cellular automata.

Let $\mu\in\Meas(\AM;\sigma)$.  For any $H\in\Natur$, we
  say $(\AM,\mu;\sigma)$ is {\dfn $(\shift{},H)$-mixing} if, for any
Borel  measurable $\bB_0,\bB_1,\ldots,\bB_H\subseteq \AM$.
\[
\lim_{\mbox{\scriptsize $\begin{array}{c}
\fm_0,\fm_1,\ldots,\fm_H\in\dM \\
|\fm_h-\fm_i|\goto\oo  \\ \forall \ h\neq i\in\CC{0...H}
\end{array}$}} 
\ \mu\lb[\Intsct_{h=0}^H \shift{-\fm_h}(\bB_h)\rb]
\quad=\quad
 \prod_{h=0}^H \mu[\bB_h].
\]
If $\dH\subset\dM$ is a finite subset,
then $\mu$ is {\dfn $\dH$-mixing} if, for any finite subset $\dB\subset\dM$
and any $\dH$-indexed
collection of $\dB$-words $\{\bb_\fh\}_{\fh\in\dH}\subset\sA^\dB$,
with cylinder sets $\bB_\fh:=\cyl{\bb_\fh}\subset\AM$ for all $\fh\in\dH$,
we have
\beqn
\label{H.mixing}
  \lim_{n\goto\oo} \ \mu\lb[\Intsct_{\fh\in\dH} \shift{-n\fh}(\bB_\fh)\rb]
\quad=\quad
\prod_{\fh\in\dH} \ \mu\lb[\bB_\fh\rb].
\eeqn
 For example, if $|\dH|=H$, then any $(\shift{},H)$-mixing measure
is $\dH$-mixing.  In particular, any nontrivial Bernoulli measure in $\AM$ is
$\dH$-mixing, and any mixing Markov measure in $\AZ$ or $\AN$ is
$\dH$-mixing.

  If $\bS\subseteq\AM$ is a subshift, then $\bS$ is {\dfn
topologically $\dH$-mixing} if, for any finite $\dB\subset\dM$ and
$\dH$-indexed collection of $\dB$-words
$\{\bb_\fh\}_{\fh\in\dH}\subseteq\bS_\dB$,
with $\bB_\fh:=\cyl{\bb_\fh}$ as above, there is some
$N\in\Natur$ such that, for all $n>N$, we have $\D \Intsct_{\fh\in\dH}
\shift{-n\fh}(\bB_\fh) \ \neq \ \emptyset$.  For example, if
$\mu$ is a $\dH$-mixing measure, then $\bS=\supp{\mu}$ is a
topologically $\dH$-mixing subshift.  In particular, any irreducible
Markov subshift of $\AZ$ or $\AN$ is topologically $\dH$-mixing.

\ignore{If $\bG\subset\AM$ is a subgroup shift, and $\bS\subset\AM$ is any 
other subshift, then it is easy to see that $\bS+\bG$ is also a subshift,
and is a disjoint union of cosets of $\bG$.}

  A {\dfn coset shift} is a subshift $\bC$ which
is a coset of some submodule shift $\bS\subset\AM$.  For example, for any
$c\in\sA$, let $c^\dM\in\AM$ denote the constant configuration equal
to $c$ everywhere.  Then $c^\dM+\bS$ is a coset shift.  If
$\bC\subseteq\AM$ is any subshift, and
$\bC-\bC:=\set{\bc-\bc'}{\bc,\bc'\in\bC}$, then it is easy to see that
\beqn
\label{coset.condition}
\statement{$\bC$ is a coset shift} \iff 
\statement{$\bC-\bC$ is a submodule shift}.
\eeqn
Let $\bc\in\AM$, and suppose $\bC:=\bc+\bS$ is a coset shift.
If $\tau^\bc:\AM\into\AM$ is the translation map
$\tau^\bc(\ba):=\ba+\bc$, then the {\dfn Haar measure} of $\bC$ is defined
then $\eta_\bC:=\tau^\bc(\eta_\bS)$, where $\eta_\bS$ is the Haar measure
of $\bS$. (This definition is independent of the choice of $\bc\in\bC$.)

 If $0\neq \varphi\in\sR$, then an $\sR$-module $\sA$ is {\dfn
$\varphi$-torsion-free} if $\varphi\,a\,\neq\,0$
 for all $a\in\sA\setminus\{0\}$.
Say $\varphi$ is a {\dfn unit} if it has a
multiplicative inverse in $\sR$.  This means that the function $\sA\ni
a\mapsto \varphi a\in\sA$ is a group automorphism of $(\sA,+)$.  If
$\sR=\Zahlmod{m}$ [Example \ref{X:LCA}(b)], then $\varphi$ is a
unit if and only if $\varphi$ is coprime to $m$.  If $\sR$ is a field
[Example \ref{X:LCA}(c)], then every nonzero element of $\sR$ is a
unit.  If $\varphi$ is a unit, then every
$\sR$-module is $\varphi$-torsion free.   

Suppose $\Phi$ has local rule (\ref{linear.CA.local.rule})
and $\sR$ has characteristic $p$.
For all $j\in\Natur$, let $\sR_j$ be the subring of $\sR$
generated by $\{\varphi_\fh^{p^j}\}_{\fh\in\dH}$.  This yields
a descending chain $\sR\supseteq\sR_1\supseteq\sR_2\supseteq\cdots$
of finite rings 
(because $\sR$ is finite), so there is some $J\in\Natur$
such that
$\sR_J=\sR_{J+1}=\cdots = \Intsct_{j=1}^\oo \sR_j$.  Let $\sR_\Phi:=\sR_J$.
For any $k\in\Natur$, 
let  $\bvphi_k:=(\sum_{\fh\in\dH} \varphi^{p^k}_\fh)-1$;
this defines a sequence $\{\bvphi_k\}_{k=1}^\oo\subseteq\sR$
which is eventually periodic (because $\sR$ is finite);
thus, there is a nonempty set $\sF_\Phi\subseteq\sR$ 
of elements which appear infinitely often in
 $\{\bvphi_k\}_{k=1}^\oo$.  We now come to our main result.

\Theorem{\label{thm:mixing.measure.rigidity}}
{
 Let $\sR$ be a finite commutative ring of prime characteristic $p$,
let $\sA$ be a finite $\sR$-module, and let $\Phi\in\RLCA[\AM]$
have  local rule {\rm(\ref{linear.CA.local.rule})}, with $|\dH|\geq 2$.
\bthmlist
\item Suppose  $\varphi_\fh$ is a unit for some $\fh\in\dH$.
If $\bC$ is a topologically $\dH$-mixing, $\Phi$-invariant subshift of $\AM$,
then $\bC$ is a coset shift of some $\sR_\Phi$-submodule shift $\bS$.

\item If there exists $\bvphi\in\sF_\Phi$ and some finite $\dB\subset\dM$ 
such that $\sA^\dB/\bS_\dB$ is $\bvphi$-torsion free
{\rm (e.g. if $\bvphi$ is a unit)}, and $\bC$ is as in {\rm(a)},
then actually $\bC=\bS$.

\item  Suppose $\varphi_\fh$ is a unit for every $\fh\in\dH$.
If $\mu$ is a $(\Phi,\shift{})$-invariant, $\dH$-mixing
measure on $\AM$, then $\mu$ is the Haar measure of a $\Phi$-invariant 
coset shift $\bC$ of some $\sR_\Phi$-submodule shift $\bS$.

If the hypothesis of {\rm(b)} holds, then $\mu$ is the Haar measure on $\bS$.
\ethmlist
}

\example{\label{X:nontrivial.invariant.coset.shift}
Let $\dM:=\Zahl\x\Natur$, let $\sR=\sA:=\Zahlmod{2}$, and let
$\dB:=\{(0,-1);\,(0,0);\,(0,1);\,(1,0)\}$; then
$\bS:=\set{\ba\in\AM}{\sum_{\fb\in\dB} a_{\fb+\fm}=0, \ \forall\, \fm\in\dM}$
is a submodule shift.  To visualize $\eta_\bS$, note that
any  $\bs\in\bS$ is entirely determined by its `zeroth row'
$\bs_{\Zahl\x\{0\}}$;  this yields a bijection $\Psi:\AZ\into\bS$,
and $\eta_\bS=\Psi(\eta)$, where $\eta$ is the uniform Bernoulli measure
on $\AZ$.  

 If $\dH:=\{(0,0);\,(0,1);\,(1,0)\}$,  then $\bS$ and $\eta_\bS$ are
$\dH$-mixing.
Let $\Phi$ be the LCA with
local rule $\phi(\ba_\dH)=\sum_{\fh\in\dH} a_\fh$
(i.e. $\varphi_\fh=1$ for all $\fh\in\dH$.).   Then $\Phi(\bS)=\bS$,
so $\Phi(\eta_\bS)=\eta_\bS$, by Proposition \ref{subgroup.haar.invariant}.
However, $\bvphi_k = |\dH|+1 = 4
\equiv 0 \pmod{2}$ for all $k\in\Natur$, so the `torsion-free'
condition of
Theorem \ref{thm:mixing.measure.rigidity}(b) is never satisfied;
thus, nontrivial $\Phi$-invariant coset shifts might exist.  Indeed,
let $\bc\in\AM$ be the `checkerboard'
configuration defined by $c_{m,n}:=(m+n)\mod{2}$, for all $(m,n)\in\dM$;
then $\bc\not\in\bS$, and $\bC:=\bc+\bS$ is a nontrivial coset shift
of $\bS$. 
Furthermore, $\Phi(\bc)=\bc$, so $\Phi(\bC)=\bC$; thus,
$\bC$ is an $\dH$-mixing, $\Phi$-invariant 
coset shift, as in Theorem \ref{thm:mixing.measure.rigidity}(a),
while $\eta_\bC$ is an $\dH$-mixing, $\Phi$-invariant 
measure, as in Theorem \ref{thm:mixing.measure.rigidity}(c).}

\Corollary{\label{rigid.cor.1}}
{
 Let $\sA=\sR:=\Zahlmod{p}$, where $p$ is prime.  Let $\dM:=\ZD\x\NE$,
and let $\Phi\in\ZpLCA$ have a neighbourhood of cardinality $H\geq 2$.
Let $\mu\in\Meas(\AM;\Phi,\shift{})$ be $(\sigma,H)$-mixing.
Suppose that either
\bdesc
  \item[{[i]}] $\DE=1$;
\qquad or
\qquad
{\bf[ii]} \ $h(\mu,\shift{})>0$ \ {\rm(and $\DE\geq1$)}.
\edesc
Then  $\mu$ is the uniform Bernoulli measure on $\AM$.
}
\bthmprf Every nonzero element of the field $\Zahlmod{p}$ is a unit,
so we can use Theorem
\ref{thm:mixing.measure.rigidity}(c).  Case {\bf[i]} is because
$(\Zahlmod{p})^\Zahl$ and $(\Zahlmod{p})^\Natur$ have no proper infinite
subgroup shifts (because if $\bS$ was such a
subgroup shift, then $\set{a\in\Zahlmod{p}}{[a,0]\in\bS_{\{0,1\}}}$
would be a proper nontrivial subgroup of $\Zahlmod{p}$, which is impossible).   Case {\bf[ii]} is because
$(\Zahlmod{p})^{\ZD\x\NE}$ has no proper subgroup shifts of nonzero
entropy \cite[first paragraph of \S25, p.228]{Schmidt95}.  \ethmprf

 Theorem \ref{thm:mixing.measure.rigidity}(c)
 is somewhat similar to \cite{Schmidt95b} or \cite[Corollary 29.5,
 p.289]{Schmidt95}, which characterizes the $\shift{}$-invariant
 measures of a subgroup shift $\bG\subseteq\AM$.  However, Schmidt
 requires $\bG$ to be `almost minimal' (i.e.  to have no infinite
 $\shift{}$-invariant subgroups), whereas we do not.

\breath

To prove Theorem \ref{thm:mixing.measure.rigidity}, 
we use tools from number theory and harmonic analysis.
If $\dM:=\ZD\x\NE$, then any $\sR$-linear CA on $\AM$ can be
written as a `Laurent polynomial of shifts' with $\sR$-coefficients.
That is, if $\Phi$ has local
rule (\ref{linear.CA.local.rule}), then for any $\ba\in\AM$,
\beqn
\label{shift.polynomial}
  \Phi(\ba)\quad:=\quad  \sum_{\fh\in\Nh} \varphi_\fh \ \shift{\fh}(\ba)
\quad\mbox{(where we add configurations componentwise).}
\eeqn
We indicate this by writing ``$\Phi=F(\shift{})$'', where 
$F\in \sR[x_1^{\pm1},\ldots,x_D^{\pm1};y_1,\ldots,y_E]$ is the 
$(\DE)$-variable Laurent polynomial defined:
\[
  F(x_1,\ldots,x_D;y_1,\ldots,y_E)\quad:=\quad 
 \sum_{(h_1,\ldots,h_D;h'_1,\ldots,h'_E)\in\Nh} \varphi_\fh \
x_1^{h_1}\ldots x_D^{h_D} y_1^{h'_1}\ldots y_E^{h'_E}.
\]
If $F$ and $G$ are two such polynomials,
and $\Phi=F(\shift{})$ while $\Gam=G(\shift{})$, then $\Phi\circ\Gam =
(F\cdot G)(\shift{})$, where $F\cdot G$ is the product of  
$F$ and $G$ in the polynomial ring
$\sR[x_1^{\pm1},\ldots,x_D^{\pm1};y_1,\ldots,y_E]$.  In particular, this
means that $\Phi^t = F^t(\shift{})$ for all $t\in\Natur$.  Thus,
iterating an $\sR$-LCA is equivalent to computing the powers of a
polynomial.  If $\sR$ is commutative, we can do this
with the Binomial Theorem, and if
$\sR$ has characteristic $p$, then we can compute the binomial coefficients
modulo $p$.  In particular, if $p$ is prime, and 
$\Phi$ has polynomial representation
(\ref{shift.polynomial}), then for any $k\in\Natur$, Fermat's Little
Theorem implies:
\beqn
\label{shift.polynomial2}
  \Phi^{p^k}\quad=\quad  \sum_{\fh\in\Nh} \varphi_\fh^{p^k} \shift{p^k\fh}.
\eeqn

\Lemma{\label{subgroupshift.lemma}}
{
  Let $\sR$ be a finite ring generated by a set $r_1,\ldots,r_H$,
at least one of which is a unit.
Let $\sA$ be an $\sR$-module and let
$\bS\subseteq\AM$ be a subshift containing $\bzero$.  The following are equivalent:
\bthmlist
\item $\bS$ is an $\sR$-submodule shift.

\item For any finite $\dB\subset\dM$, \ $\bS_\dB$ is an $\sR$-submodule
of $\sA^\dB$.

\item For any finite $\dB\subset\dM$ and any
$\bs_1,\ldots,\bs_H\in\bS_\dB$, we have
$\D \sum_{h=1}^H r_h \bs_h \ \in \ \bS_\dB$.
\ethmlist
}
\bthmprf
 Clearly $\mathbf{(a)}\Leftrightarrow\mathbf{(b)}\Rightarrow\mathbf{(c)}$. 
To see that
$\mathbf{(c)}\Rightarrow\mathbf{(b)}$, 
let $\tlsR$ be the set of $\Natur$-linear combinations of products
of $\{r_h\}_{h=1}^H$.  Then {\bf(c)} implies that 
$\tlsR\bS_\dB\subseteq\bS_\dB$ (because $\bzero\in\bS$).  Thus
it suffices to show that $\tlsR=\sR$.

  Every element of $\sR$ is a
$\Zahl$-linear combination of products of $\{r_h\}_{h=1}^H$; hence we
need only show that $-1_\sR\in\tlsR$.  Some $r_h$ is a unit, and
$\sR$ is finite, so there exists $n\in\Natur$ with
$1_\sR=r_h^n\in\tlsR$.  But $\sR$ has  characteristic $c<\oo$, so
$-1_\sR=(c-1)1_\sR \in\tlsR$. \ethmprf

\bthmprf[Proof of Theorem \ref{thm:mixing.measure.rigidity}{\rm(a)}.]
Let $\bC\subset\AM$ be a topologically $\dH$-mixing,
$\Phi$-invariant subshift.   Let $\bS:=\bC-\bC$; then 
Fact (\ref{coset.condition}) says that
it suffices to show that $\bS$ is an $\sR_\Phi$-submodule shift.
Now, $\bS$ is a subshift, $\bzero\in\bS$, and for all $k\in\Natur$,
some element of $\{\varphi^{p^k}_\fh\}_{\fh\in\dH}$ is a unit, so we
can use Lemma \ref{subgroupshift.lemma}.

Let $\dB\subset\dM$ be finite.  Now, $\bS$ is $\dH$-mixing
(because $\bC$ is), so there exists $N\in\Natur$ such that, for any $n>N$,
and any $\dH$-indexed collection $\{\bb_\fh\}_{\fh\in\dH}\subseteq\bS_\dB$, 
there exists $\bs\in\bS$ with
\beqn
\label{thm:mixing.measure.rigidity.i.e1}
\bs_{\dB+n\fh} \ = \ \bb_\fh, \ \ \mbox{for all $\fh\in\dH$}.  
\eeqn
  Make $k$ large enough that $p^k>N$, and such that
$\sR_\Phi$ is generated by $\{\varphi_\fh^{p^k}\}_{\fh\in\dH}$.
Let $\{\bb_\fh\}_{\fh\in\dH}\subseteq\bS_\dB$ be arbitrary,
and find  
$\bs\in\bS$ satisfying eqn.(\ref{thm:mixing.measure.rigidity.i.e1})
for $n:=p^k$.
Then
\[
  \Phi^{p^k}(\bs)_\dB\quad\eeequals{(*)}\quad  \sum_{\fh\in\Nh} \varphi_\fh^{p^k} \shift{p^k\fh}(\bs)_\dB
\quad=\quad  \sum_{\fh\in\Nh} \varphi_\fh^{p^k}\,\bs_{\dB+p^k\fh}
\quad\eeequals{(\dagger)}\quad
\sum_{\fh\in\Nh} \varphi_\fh^{p^k} \bb_\fh.
\]
[$(*)$ is by eqn.(\ref{shift.polynomial2}) and
$(\dagger)$ is by eqn.(\ref{thm:mixing.measure.rigidity.i.e1}).]
But $[\Phi^{p^k}(\bs)]_\dB \in \bS_\dB$, because 
$\Phi(\bS)\subseteq\bS$, because 
$\Phi(\bC)\subseteq\bC$.  This verifies condition (c) of
Lemma \ref{subgroupshift.lemma}
for any finite $\dB\subset\dM$ and 
$\{\bb_\fh\}_{\fh\in\dH}\subseteq\bS_\dB$.
\ethmprf

\Lemma{\label{inv.coset.shift}}
{
  Let $\sR$, $\sA$, $\Phi$, and $\sF_\Phi$
be as in {\rm Theorem \ref{thm:mixing.measure.rigidity}}, and let $\bS\subset\AM$ be a submodule shift.   If $\bC$ is any
$\dH$-mixing, $\Phi$-invariant coset shift of $\bS$,
 then $\bvphi\,\bC\ \subseteq \ \bS$ for all $\bvphi\in\sF_\Phi$.
}
\bthmprf
Let $\bvphi\in\sF_\Phi$ and let $\ba\in\bC$.  To show that
$\bvphi\,\ba\ \in \ \bS$, it suffices to show, 
for any finite $\dB\subset\dM$, that
$\bvphi\,\ba_\dB \ \in\ \bS_\dB$.

There exist arbitrarily large $k\in\Natur$ with $\bvphi_k=\bvphi$.
But if $k$ is large enough, then there exists $\bc\in\bC$ with
$\bc_{p^k\fh+\dB}=\ba_\dB$ for all $\fh\in\dH$ (because $\bC$ is
$\dH$-mixing).  Thus, $\Phi^{p^k}(\bc)_\dB= (\sum_{\fh\in\dH}
\varphi_\fh^{p^k}) \, \ba_\dB$, by eqn.(\ref{shift.polynomial2}).  Thus,
$ \bvphi_k \, \ba_\dB \ = \ \Phi^{p^k}(\bc)_\dB-\bc_\dB$.  But
$\Phi^{p^k}(\bc)_\dB-\bc_\dB \ = \  (\Phi^{p^k}(\bc)-\bc)_\dB \ \in \ \bS_\dB$
because $\Phi^{p^k}(\bc)-\bc\ \iiin{(*)}\ \bC-\bC
\ \eeequals{(\dagger)} \ \bS$, where $(*)$ is because
$\Phi^{p^k}(\bC)\subseteq\bC$, and $(\dagger)$ is by Fact
(\ref{coset.condition}).  \ethmprf

\bthmprf[Proof of Theorem \ref{thm:mixing.measure.rigidity}{\rm(b)}.]
(by contradiction) \ 
If $\AB/\bS_\dB$ is $\bvphi$-torsion free, then
$\sA^{\dB'}/\bS_{\dB'}$ is also $\bvphi$-torsion free,
 for any $\dB'\supseteq\dB$.
If $\bC\neq\bS$, and $\dB$ is large enough,
then $\bC_\dB$ is a nontrivial coset of the submodule $\bS_\dB$ 
in $\AB$; thus, $\bC_\dB$ is nontrivial as an element of
the quotient module $\sA^\dB/\bS_\dB$.
But Lemma \ref{inv.coset.shift}
implies that $\bvphi\bC_\dB\subseteq \bS_\dB$, so
$\bvphi$ annihilates $\bC_\dB$ in $\sA^\dB/\bS_\dB$, so
$\sA^\dB/\bS_\dB$ has nontrivial $\bvphi$-torsion,
which is a contradiction. By contradiction, $\bC=\bS$.
\ethmprf

  Let $\dT:=\set{c\in\Cplx}{|c|=1}$ be the
unit circle group.  If $(\bG,+)$ is a compact abelian group
(e.g. $\bG:=\AM$ where $(\sA,+)$ is a finite abelian group),
then a {\dfn character} on $\bG$ is a continuous group
homomorphism $\chi:(\bG,+)\into(\dT,\cdot)$.  Let
$\hbG$ denote the group of characters of $\bG$.
If $\mu\in\Meas(\bG)$, then $\mu$ is uniquely identified by its
{\dfn Fourier coefficients}
\[
 \hmu[\chi]
\quad := \quad \int_{\bG} \chi \ d \mu,\qquad \mbox{for all $\chi\in\hbG$}.
\]
For example, if  $\eta_\bG$ is the Haar measure on $\bG$, 
and $\chr{}\in\hbG$ is the trivial character, then it is 
easy to verify:

\Lemma{\label{haar.lemma0}}
{
$\eta_\bG$ is the unique Borel measure on $\bG$
such that $\heta_\bG[\chr{}]=1$ and $\heta_\bG[\chi]=0$ for all
other $\chi\in\hbG$.\qed}

  More generally, we have the following:

\Lemma{\label{Haar.lemma}}
{
  Let $(\bG,+)$ be a compact abelian group, and let $\mu
\in\Meas(\bG)$. Then
\\
$\statement{$\mu=\eta_\bS$ for some closed subgroup $\bS\subseteq\bG$}$
$\iff
\statement{$\forall \ \chi\in\hbG$, either $\hmu[\chi]=0$ or $\hmu[\chi]=1$}$.

}
\bthmprf
``$\implies$'' \ Suppose $\mu=\eta_\bS$.  If $\chi\in\hbG$, then
$\chi\restr{\bS}\in\hbS$.  Thus, Lemma \ref{haar.lemma0} implies that
\[
\hmu[\chi]\quad=\quad\hmu\lb[\chi\restr{\bS}\rb]\quad=\quad
\choice{1 &\If& \chi\restr{\bS} \equiv \chr{\bS}; \\
        0 &&\mbox{otherwise}.}
\]
``$\Leftarrow$'' \ If $\gam\in\hbG$ and
$\hmu[\gam]=1$, then $\supp{\mu}\subseteq\ker(\gam)$.
If $\bS:=\Intsct\lb\{\ker(\gam)\,;\, \gam\in\hbG,\, \hmu[\gam]=1\rb\}$,
then $\bS$ is a closed subgroup of $\bG$, and $\supp{\mu}\subseteq\bS$.
We claim $\mu=\eta_\bS$.
If $\chi\in\hbS$, then $\chi=\gam\restr{\bS}$ for some
$\gam\in\hbG$ (this follows from the Pontrjagin Duality Theorem;
see e.g. Fact (6), \S0.7, p.13 of \cite{Walters}).
If $\chi\neq\chr{\bS}$, then $\bS\not\subset\ker(\gam)$,
so $\hmu[\gam]\neq 1$, so $\hmu[\gam]=0$ (by hypothesis), so $\hmu[\chi]=0$.
Thus, Lemma \ref{haar.lemma0} implies that $\mu=\eta_\bS$.
\ethmprf

  Let $(\sA,+)$ be a finite abelian group, so that $(\AM,+)$ is
compact abelian.  For any $\bchi\in\hAM$, there is a unique finite
subset $\dB\subset\dM$ and  unique  nontrivial
$\chi_\fb\in\hsA$ for all $\fb\in\dB$ such that
\beqn
\label{character.eqn}
  \bchi(\ba)\quad=\quad\prod_{\fb\in\dB} \chi_\fk(a_\fb),
\qquad\forall \ \ba\in\AM.
\eeqn  
  We say that $\bchi$ is {\dfn based} on $\dB$.  If 
$\{\bchi^\fh\}_{\fh\in\dH}\subset\hAM$
is an $\dH$-indexed collection of characters based on $\dB$, 
and $\mu$ is $\dH$-mixing, then
equations (\ref{H.mixing}) and (\ref{character.eqn}) together imply
\beqn
\label{H.mixing2}
  \lim_{n\goto\oo} \ \hmu\lb[\prod_{\fh\in\dH} \bchi^\fh\circ\shift{-n\fh}\rb]
\quad=\quad
\prod_{\fh\in\dH} \ \hmu\lb[\bchi^\fh\rb].
\eeqn

\bthmprf[Proof of Theorem \ref{thm:mixing.measure.rigidity}{\rm(c)}.]
If  $\bC=\supp{\mu}$, then
Theorem \ref{thm:mixing.measure.rigidity}{\rm(a)}
says $\bC$ is an $\sR_\Phi$-coset shift ---i.e.
$\bC=\bc+\bS$, where $\bS$ is an $\sR_\Phi$-submodule shift.
Let $\nu:=\tau^{-\bc}(\mu)$; then $\supp{\nu}=\bS$.
To show that $\mu=\eta_\bC$,
we must show that $\nu=\eta_{\bS}$; we will do this with
Lemma \ref{Haar.lemma}.  For any $\bchi\in\hAM$, note that
\beqn
\lb|\hmu[\bchi]\rb|
\quad\eeequals{(*)}\quad
\lb|\bchi(\bc)\cdot\hnu[\bchi]\rb|
\quad=\quad
\lb|\bchi(\bc)\rb|\cdot\lb|\hnu[\bchi]\rb|
\quad\eeequals{(\dagger)}\quad
\lb|\hnu[\bchi]\rb|.
\label{thm:mixing.measure.rigidity.e0}
\eeqn
[Here, $(*)$ is because $\mu=\tau^\bc(\nu)$ and  $\bchi$ is a homomorphism, 
while $(\dagger)$ is because $\bchi[\bc]\in\dT$.] \
Now, let $\bchi\in\hAM$, and suppose
$\hnu[\bchi]\neq 1$; we must show that $\hnu[\bchi]= 0$.
Let $k$ be large enough that $\sR_\Phi$ is generated by 
$\{\varphi^{p^k}_\fh\}_{\fh\in\dH}$.  Let $\Gam$ be the multiplicative
group generated by $\{\varphi^{p^k}_\fh\}_{\fh\in\dH}$
(which are units, by hypothesis).

\Claim{\label{thm:mixing.measure.rigidity.C1}  $\Gam\subset\sR_\Phi$.}
\bclaimprf
By definition, $\sR_\Phi$ contains all products of positive powers of 
$\{\varphi^{p^k}_\fh\}_{\fh\in\dH}$; we must show it contains the
negative powers as well.  For all $\fh\in\dH$, the unit
 $\varphi^{p^k}_\fh$ has finite multiplicative order,
because $\sR_\Phi$ is finite; thus, there is some $n\in\Natur$
such that $(\varphi^{p^k}_\fh)^{-1}=(\varphi^{p^k}_\fh)^n\in\sR_\Phi$.
Thus, $\sR_\Phi$ contains all products of integer powers 
(positive or negative) of 
$\{\varphi^{p^k}_\fh\}_{\fh\in\dH}$; hence $\Gam\subset\sR_\Phi$.
\eclaimprf
If $\gam\in\Gam$, then $\bchi\circ\gam\in\hAM$ 
(it is a composition of two homomorphisms).  Define
\beqn
\label{M.defn}
  M\quad:=\quad \max_{\gam\in\Gam} \ \lb|\hnu[\bchi\circ\gam]\rb|.
\eeqn
\ignore{
\[
  M\quad:=\quad \max\set{\lb|\hnu[\bchi\circ\varphi_\fh^n ]\rb|}
{\fh\in\dH \And n\in\CO{0...p}}.
\]}
Thus, $\lb|\hnu[\bchi]\rb|\leq M$.  We will show that $M=0$.

\Claim{\label{thm:mixing.measure.rigidity.C2} $M<1$.}
\bclaimprf
 (by contradiction). \  $\Gam$ is finite, so if
$M=1$, then there is some
$\gam\in\Gam$  such that
 $\lb|\hnu[\bchi\circ\gam ]\rb|=1$.  Let $t:=\hnu[\bchi\circ\gam]$;
then $t\in\dT$, and $\supp{\nu}\subseteq (\bchi\circ\gam)^{-1}\{t\}$.
But $\bzero\in\supp{\nu}$ (because $\supp{\nu}=\bS$ is a submodule),
so $t=1$; hence $\hnu[\bchi\circ\gam ]=1$.
Thus, $\nu\lb[\ker(\bchi\circ\gam)\rb]=1$, 
hence $\bS\subseteq\ker(\bchi\circ\gam)$,
hence $\gam(\bS)\subseteq\ker(\bchi)$.

\quad Now $\gam^{-1}\in\Gam\subset\sR_\Phi$ (by
 Claim \ref{thm:mixing.measure.rigidity.C1}), and $\bS$
 is an $\sR_\Phi$-module (by Theorem \ref{thm:mixing.measure.rigidity}(a)), so 
$\gam^{-1}(\bS)\subseteq\bS$, so
$\bS\subseteq \gam(\bS) \, \subseteeeq{(*)} \,
\ker(\bchi)$, where $(*)$ is by the previous paragraph.  But then
$\hnu[\bchi]=1$, which contradicts the definition of $\bchi$.

By contradiction, we must have $M<1$.
\eclaimprf

  By replacing $\bchi$ with $\bchi\circ\gam$ for some $\gam\in\Gam$
if necessary,
we can assume without loss of generality that  $\lb|\hnu[\bchi]\rb|= M$. 
Then $\lb|\hmu[\bchi]\rb|= M$ also, 
by eqn.(\ref{thm:mixing.measure.rigidity.e0}). 
Thus it suffices to evaluate $\hmu[\bchi]$.
But if $\Phi$ has polynomial representation (\ref{shift.polynomial}),
then for any $k\in\Natur$, we have
\begin{eqnarray}
\nonumber
\hmu[\bchi]&:=&
\int_{\AM} \bchi \dmu
\quad\eeequals{(a)}\quad \int_{\AM} \bchi \ d(\Phi^{(p^k)}\mu)
\quad\eeequals{(b)}\quad \int_{\AM} \bchi\circ \Phi^{(p^k)} \dmu
\nonumber
\\ &\eeequals{(c)} &
 \int_{\AM}  \bchi\circ \lb(\sum_{\fh\in\dH}\varphi_\fh^{(p^k\fh)}\shift{(p^k\fh)}\rb) \dmu
\quad\eeequals{(h)}\quad
 \int_{\AM}  \prod_{\fh\in\dH} \lb( \bchi\circ\varphi_\fh^{(p^k\fh)}\circ\shift{(p^k\fh)}\rb) \dmu\nonumber
\\&=&
\hmu\lb[ \prod_{\fh\in\dH} \bchi\circ\varphi_\fh^{(p^k\fh)}\circ\shift{(p^k\fh)} \rb].
\label{thm:mixing.measure.rigidity.e1}
\end{eqnarray}
[Here, (a) is because $\mu$ is $\Phi$-invariant; \  (b) is a change of variables; \  (c) is by eqn.(\ref{shift.polynomial2});  and (h) is because
$\bchi$ is a homomorphism.]  \
Thus,
\begin{eqnarray}\hmu[\bchi] &\eeequals{(\ref{thm:mixing.measure.rigidity.e1})}&
\lim_{k\goto\oo} \ 
\hmu\lb[ \prod_{\fh\in\dH} \bchi\circ\varphi_\fh^{(p^k\fh)}\circ\shift{(p^k\fh)} \rb]
\quad\leeeq{(*)} \quad
\limsup_{k\goto\oo} \ 
\prod_{\fh\in\dH} \hmu\lb[\bchi\circ\varphi_\fh^{(p^k\fh)}\rb],
\label{thm:mixing.measure.rigidity.e2}
\end{eqnarray}
where $(*)$ is by eqn.(\ref{H.mixing2}), because $\mu$ is $\dH$-mixing
and $\Gam$ is finite.  Thus,
\[
M \quad=\quad
\lb|\hmu[\bchi]\rb|
\quad\leeeq{(\ref{thm:mixing.measure.rigidity.e2})}\quad
\limsup_{k\goto\oo} \ \prod_{\fh\in\dH} \lb| \hmu\lb[ \bchi\circ\varphi_\fh^{(p^k\fh)} \rb]\rb|
\quad \leeeq{(*)}  \quad 
\lim_{k\goto\oo} \ \prod_{\fh\in\dH} M
\quad=\quad M^H.
\]
[Here, $(*)$ is by 
equations (\ref{thm:mixing.measure.rigidity.e0}) and (\ref{M.defn}).] \ 
But $H\geq 2$, and Claim \ref{thm:mixing.measure.rigidity.C2}
says $M<1$.  Thus, $M=0$.  Thus, $\lb|\hnu[\bchi]\rb|=0$.  

This
holds for any $\bchi\in\hAM$ with $\hnu[\bchi]\neq 1$, so 
Lemma \ref{Haar.lemma} says $\nu=\eta_\bS$; hence $\mu=\eta_\bC$.
\ethmprf

\breath

{\em Constructing coset shifts.}  \
To enumerate all $\dH$-mixing, $\Phi$-invariant measures, 
Theorem \ref{thm:mixing.measure.rigidity}(c) says it suffices to
enumerate all $\dH$-mixing, $\Phi$-invariant coset shifts.
But if $\bS\subset\AM$ is a submodule shift, not
every coset of $\bS$ is a coset {\em shift}.  Indeed,
let $\bc\in\AM$, and
for all $\fm\in\dM$, let $\bb^\fm:=\shift{\fm}(\bc)-\bc$.
Then it is easy to check
\beqn
\label{coset.shift.condition}
\statement{The coset $\bc+\bS$ is a subshift}
\iff
\statement{$\bb^\fm\in\bS$ for all $\fm\in\dM$}.
\eeqn
We call $\bB:=\{\bb^\fm\}_{\fm\in\dM}$
the {\dfn coboundary} of $\bc$.  More generally,
an $\AM$-valued {\dfn cocycle} is any $\dM$-indexed collection  
$\{\bb^\fm\}_{\fm\in\dM}\subseteq\AM$ such that
$\bb^{\fm+\fn} 
\ =\  \shift{\fn}(\bb^\fm) + \bb^\fn$, for all
$\fn,\fm\in\dM$.

\Lemma{\label{cocycle.lemma}}
{
\bthmlist
  \item If $\bc\in\AM$, then its coboundary $\bB$ is a cocycle.

  \item $\bc$ is entirely determined by $c_0$ and $\bB$, because
 $c_\fm= c_0 + b^\fm_0$, for all $\fm\in\dM$.

  \item Let $\bB:=\{\bb^\fm\}_{\fm\in\dM}\subset\bS$ be any cocycle,
and let $c_0\in\sA$ be arbitrary.  Define $\bc:=[c_\fm]_{\fm\in\dM}$
according to part {\rm(b)}.  Then $\bB$ is the coboundary of $\bc$,
\ethmlist
}
\bthmprf {\bf(a)} is straightforward.  {\bf(b)} is
because $c_\fm=(\shift{\fm}(\bc))_0 = (\bc+\bb^\fm)_0 = c_0 +
b^\fm_0$.  For {\bf(c)}, let $\fm,\fn\in\dM$.  Then
$[\shift{\fm}(\bc)-\bc]_\fn  \ = \  c_{\fm+\fn} - c_\fn
\ := \ (c_0 + b^{\fm+\fn}_0) - (c_0 + b_0^\fn)
\ = \ b^{\fm+\fn}_0  - b_0^\fn
\ = \ (\bb^{\fm+\fn}  - \bb^\fn)_0
\ \eeequals{(*)} \ \shift{\fn}(\bb^\fm)_0
\ = \ b^\fm_\fn$, where $(*)$ is the cocycle property.
This holds for all $\fn\in\dM$, so $\shift{\fm}(\bc)-\bc=\bb^\fm$.
This holds for all $\fm\in\dM$, so $\bB$ is the coboundary of $\bc$.
\ethmprf

Lemma \ref{cocycle.lemma} and Fact (\ref{coset.shift.condition}) imply
that, to construct a coset shift of $\bS$, it suffices to construct
an $\bS$-valued cocycle. In particular, if $a\in\sA$ such that
$a^\dM\in\bS$, then we can get one $\bS$-valued cocycle by defining
$\bb^\fm:=(m_1+\cdots+m_D) a^\dM$, for all $\fm=(m_1,\ldots,m_D)\in\dM$;
in Lemma \ref{cocycle.lemma}(c) in this case, 
$c_\fm=c_0+(m_1+\cdots+m_D) a$ for all $\fm\in\dM$.
For example, if $\bS$ is as in Example
\ref{X:nontrivial.invariant.coset.shift}, then $1^\dM\in\bS$.  In 
this case, $\bb^{(m,n)}=0^\dM$ if $m+n$ is even, and $\bb^{(m,n)}=1^\dM$ if
$m+n$ is odd.  If $c_0:=0$, and we define $\bc$ as in Lemma
\ref{cocycle.lemma}(b), then we get precisely the `checkerboard'
configuration of Example \ref{X:nontrivial.invariant.coset.shift};
thus $\bC=\bc+\bS$ is a nontrivial coset shift, as previously claimed.

\breath

{\em Extension to rings of squarefree characteristic.} \ 
  We say $m\in\Natur$ is {\dfn squarefree} if $m$ has prime
factorization $m=p_1p_2\cdots p_J$, where $p_1,\ldots,p_J$ are
distinct primes.  If
$\sR$ has {\em nonprime} characteristic, then
eqn.(\ref{shift.polynomial2}) is no longer true.  However, the next
result allows us to reduce squarefree-characteristic LCA
to the previous case of prime-characteristic LCA.

\Proposition{\label{CRT}}
{
  Let $\sR$ be a commutative ring with
characteristic $m=p_1^{s_1}p_2^{s_2}\cdots p_J^{s_J}$,
where $p_1,\ldots,p_J$ are distinct primes.  For all $j\in\CC{1...J}$,
let $q_j:= m/p_j^{s_j}$, and let $\sI_j:=\set{r\in\sR}{q_j r=0}$.
Then
\bthmlist
  \item $\sI_j$ is an ideal of $\sR$, and the quotient ring
 $\sR_j:=\sR/\sI_j$ has characteristic $p_j^{s_j}$.

  \item $\sR$ is isomorphic to the direct product
$\tlsR:=\sR_1\dirsum\sR_2\dirsum\cdots\dirsum\sR_J$, via the map
$\sR\ni r\mapsto (r^1,r^2,\cdots,r^J)\in\tlsR$, where
for all $j\in\CC{1..,J}$, we define $r^j := (r+\sI_j)\in\sR_j$.

  \item Let $\sA$ be any $\sR$-module.  For all $j\in\CC{1...J}$,
let $\sB_j:=\sI_j\sA :=\set{i a}{i\in\sI_j \And a\in\sA}$ {\rm(a submodule)}, and
let $\sA_j:=\sA/\sB_j$ be the quotient module.  Then
 $\sA\cong \tlsA:=\sA_1\dirsum\cdots\dirsum\sA_J$, via the map
$\sA\ni a\mapsto (a^1,\cdots,a^J)\in\tlsA$, where
for all $j\in\CC{1..,J}$, we define $a^j := (a+\sB_j)\in\sA_j$.
  Furthermore, $\tlsR$ acts on $\tlsA$ componentwise; that is, for
any $\tlr=(r^1,\ldots,r^J)\in\tlsR$ and $\tla=(a^1,\ldots,a^J)\in\tlsA$,
we have $\tlr\cdot \tla \ = \ (r^1a^1,\ldots,r^Ja^J)$.
\ethmlist
}
\bthmprf
{\bf(a)} is straightforward. 
 {\bf(b)} follows from the Chinese Remainder Theorem for rings 
\cite[Thm.17, p.268, \S10.3]{DummitFoote} and the following claim.

\Claim{
{\bf[i]} \ $\sI_1\intsct\sI_2\intsct\cdots\intsct\sI_J=\{0\}$,
\ and \
{\bf[ii]} \ For any $j\neq k$, \ $\sI_j+\sI_k=\sR$.
}
\bclaimprf
{\bf[i]} If $r\in\Intsct_{j=1}^J \sI_j$ then $q_jr=0$ for all $j\in\CC{1...J}$.
But $\gcd(q_1,\ldots,q_J)=1$, so Bezout's identity yields
$z_1,\ldots,z_J\in\Zahl$ such that $\sum_{j=1}^J z_j q_j=1$.
But then $r= 1\cdot r = \sum_{j=1}^J z_j q_jr = \sum_{j=1}^J z_j 0 = 0$.

{\bf[ii]}  $\sI_j+\sI_k$ is an ideal, so
it suffices to show that $\sI_j+\sI_k$ contains $1_\sR$.   Let $\tlbZ_m$ be
the subring of $\sR$ generated by $1_\sR$; then $\tlbZ_m$ is isomorphic
to $\Zahlmod{m}$.  It is easy to check that
$\tlbZ_m\intsct \sI_j = p_j^{s_j} \tlbZ_m$ and
$\tlbZ_m\intsct \sI_k = p_k^{s_k} \tlbZ_m$.  But then
$1_\sR \in p_j^{s_j} \tlbZ_m + p_k^{s_k} \tlbZ_m$,
because $1 \in p_j^{s_j} \Zahlmod{m} + p_k^{s_k} \Zahlmod{m}$
by Bezout's identity, because $\gcd( p_j^{s_j}, p_k^{s_k})=1$
(because $p_j$ and $p_k$ are distinct primes).
\eclaimprf

{\bf(c)} follows from Claim 1 and
the Chinese Remainder Theorem for modules
\cite[Ex.16-17, p.333, \S10.3]{DummitFoote}.
\ethmprf

\example{Suppose $\sA=\sR:=\Zahlmod{m}$, where $m=p_1\cdots p_J$. 
 Then for all $j\in\CC{1...J}$, \ we have
$\sI_j=p_j \Zahlmod{m}$ in Proposition \ref{CRT}(a),
so $\sA_j=\sR_j=\Zahlmod{m}/(p_j \Zahlmod{m})
\cong \Zahlmod{p_j}$ in Proposition \ref{CRT}(b,c)
which is the classic Chinese
Remainder Theorem.}

Let  $\sA\cong \D \Dirsum_{j=1}^J\sA_j$
as in Proposition \ref{CRT}(c).
Define $\Psi:\AM \into\D  \Dirsum_{j=1}^J \sA_j^\dM$
by 
\beqn
\label{psi.defn}
\Psi\lb([a_\fm]_{\fm\in\dM}\rb) \ := \
\lb([a^1_\fm]_{\fm\in\dM},[a^2_\fm]_{\fm\in\dM},\ldots,[a^J_\fm]_{\fm\in\dM}\rb),
\eeqn
where, for any $a\in\sA$, we write $a \cong (a^1,\cdots,a^J)$ as in
Proposition \ref{CRT}(c).  Then $\Psi$ is a $\shift{}$-commuting,
homeomorphic $\sR$-module isomorphism.  If $\mu\in\Meas(\AM)$, then
for all $j\in\CC{1...J}$, let $\mu_j$ be the projection of $\Psi(\mu)$ to
$\AjM$; we say that $\mu$ is a {\dfn joining} of $\mu_1,\ldots,\mu_J$.
(See e.g. \cite{delaRue} or \cite[Ch.6]{Rudolph} for more about
joinings.)

\Proposition{\label{thm:mixing.measure.rigidity2}}
{
 Let $\sR$ be a commutative ring of squarefree characteristic,
let $\sA$ be an $\sR$-module, and write $\sA\cong\Dirsum_{j=1}^J\sA_j$
as in {\rm Proposition \ref{CRT}(c)}.  Let $\Phi\in\RLCA[\AM]$
have neighbourhood $\dH$ and coefficients $\{\varphi_\fh\}_{\fh\in\dH}$,
all of which are units.
If $\mu\in\Meas(\AM;\Phi,\shift{})$ is $\dH$-mixing,
then $\mu$ is a joining of measures $\mu_1,\ldots,\mu_J$, where
for each $j\in\CC{1...J}$, $\mu_j$ is the Haar measure of some
$\Phi_j$-invariant $\sR_j$-coset shift of $\AjM$.
}
\bthmprf
  For all $\fh\in\dH$, write $\varphi_\fh=(\varphi^1_\fh,\ldots,\varphi^J_\fh)$
as in Proposition \ref{CRT}(b).    Then $\varphi^j_\fh$ is a unit in $\sR_j$
for each $j\in\CC{1...J}$.
For all $j\in\CC{1...J}$, let
$\Phi_j\in\RjLCA$ have local rule
$\phi(\ba^j_\Nh) \:=\ \sum_{\fh\in\Nh} \varphi^j_\fh a^j_\fh$,
for any $\ba^j_\Nh \in\sA_j^\Nh$.
If $\Psi$ is as in eqn.(\ref{psi.defn}), then
Proposition \ref{CRT}(c) implies that $\Psi$ is a topological
conjugacy from $(\AM,\Phi)$ to
the direct product $(\sA_1^\dM,\Phi_1)\x\cdots\x(\sA_J^\dM,\Phi_J)$.
For all $j\in\CC{1...J}$, let $\mu_j$ be the projection of $\Psi(\mu)$ to 
$\AjM$; then  $\mu_j\in\Meas(\AjM;\Phi_j,\shift{})$ and is $\dH$-mixing; 
hence, Theorem \ref{thm:mixing.measure.rigidity}(c) implies
that $\mu_j$ is the Haar measure for some  $\sR_j$-coset shift of $\AjM$.
\ethmprf

\Corollary{\label{rigid.cor.3}}
{
 Let $\sA=\Zahlmod{m}$, where $m=p_1\cdots p_J$ is squarefree.
Let $\dM:=\ZD\x\NE$, and
let $\Phi\in\ZmLCA$ have a neighbourhood of cardinality $H\geq 2$.
Let $\mu\in\Meas(\AM;\Phi,\shift{})$ be $(\sigma,H)$-mixing, and
suppose that either:
\bdesc
  \item[{[i]}] $D+E=1$;
\qquad or
\qquad
  {\bf[ii]} \ $h(\mu,\shift{}) \ >  \ 
\log_2(m)-\log_2 \lb(\min\{p_1,\ldots,p_J\}\rb)$. 
\edesc
Then $\mu$ is a joining of
the uniform Bernoulli measures on $\sA_1^\dM,\ldots,\sA_J^\dM$,
where $\sA_j:=\Zahlmod{p_j}$ for each $j\in\CC{1...J}$.
}
\bthmprf Case {\bf[i]} follows immediately from 
Proposition \ref{thm:mixing.measure.rigidity2} and Case [i] of Corollary \ref{rigid.cor.1}. 

Case {\bf[ii]}: \ 
For each $i\in\CC{1...J}$, let
$\mu_i$ be as in Proposition \ref{thm:mixing.measure.rigidity2}.
Then 
\beq
h(\mu_i) &\geeeq{(*)}&
h(\mu)-\sum_{i\neq j=1}^J h(\mu_j)
\quad\geeeq{(\dagger)}\quad
h(\mu)-\sum_{i\neq j=1}^J \log_2(p_j)
\quad\eeequals{(\ddagger)}\quad
h(\mu)-\log_2(m)+\log_2(p_i)
\\&\grt{(\diamond)}&
\log_2(m)-\log_2 \lb(\min\{p_j\}_{j=1}^J\rb)
-\log_2(m)+\log_2(p_i)
\quad\geq\quad 0.
\eeq
Here, $(*)$ is because $h(\mu)\leq\sum_{j=1}^J h(\mu_j)$, while
$(\dagger)$ is because for each $j\in\CC{1...J}$ we have
 $h(\mu_j)\leq \htop(\AjM)=\log_2|\sA_j|=
\log_2(p_j)$.  Next, $(\ddagger)$ is because $m=p_1\cdots p_J$, and
$(\diamond)$ is by the hypothesis of Case {\bf[ii]}. 

Thus, Case [ii] of Corollary \ref{rigid.cor.1} 
implies that $\mu_i$ is the uniform Bernoulli measure on $\AiM$.
This holds for all $i\in\CC{1...J}$; the result follows.
\ethmprf

\ignore{
\noindent {\em Question.}\  For all $j\in\CC{1...J}$, let
$\mu_j\in\Meas(\AjM;\Phi_j,\shift{})$ be as in Proposition
\ref{thm:mixing.measure.rigidity2}.   Aside from the product measure
$\mu_1\x\cdots\x\mu_J$, what other  joinings of
$\mu_1,\ldots,\mu_J$ exist on $\AM$?}

{\footnotesize
\bibliographystyle{alpha}
\bibliography{bibliography}
}
\end{document}